\documentclass[11pt]{article}
\usepackage[cm]{fullpage}

\usepackage{amsfonts}
\usepackage{amsmath}
\usepackage{amssymb}
\usepackage{hyperref}
\usepackage{mathrsfs}
\usepackage[all,cmtip]{xy}
\usepackage{enumitem}
\usepackage{ulem}
\usepackage{eufrak}
\usepackage[all,cmtip]{xy}

\usepackage{amsmath, amsthm, amsfonts, amssymb}
\usepackage{graphicx}
\usepackage{subcaption}

\usepackage{amsmath,amscd, float,rotating}
\usepackage{pb-diagram}

\def\endpf{\hbox{\vrule height1.5ex width.5em}}

\newcommand{\dbar}{\overline{\partial}}

\newcommand{\ddbar}{\sqrt{-1}\partial\dbar}

\newtheorem{remark}{Remark}[section]

\newtheorem{theorem}{Theorem}[section]
\newtheorem{proposition}{Proposition}[section]
\newtheorem{lemma}{Lemma}[section]
\newtheorem{example}{Example}[section]

\newtheorem{definition}{Definition}[section]

\newtheorem{question}{Question}[section]

\begin{document}

\title{\bf On the construction of a complete K\"ahler-Einstein metric with negative scalar curvature near an isolated log-canonical singularity}

\author{ Hanlong Fang$^{}$   and  Xin Fu$^{}$}
\date{}

\maketitle
\bigskip
\bigskip
\noindent
{\bf Abstract} \ In this short note we are concerned with the K\"ahler-Einstein metrics near  cone type log canonical singularities. 
By two different approaches, we construct a complete K\"ahler-Einstein metric with negative scalar curvature in a neighborhood of the cone  over 
a Calabi-Yau manifold, which provides a local model for the future study of the global K\"ahler-Einstein metrics on singular varieties.  In the first approach, we show that the singularity is uniformized by a complex ball and hence the induced metric from the Bergman metric of the ball is a desired one. In the second approach, we obtain a complete K\"ahler-Einstein metric with negative curvature by using Calabi Ansatz. At last, we show that these two metrics are indeed the same.


\section{Introduction}\label{intro}
Existence of K\"ahler-Einstein metric has been the central topic in complex geometry for decades. For the canonically polarized compact manifolds, the existence was established by \cite{A,Y1};  Calabi-Yau case was solved by \cite{Y1}. Also, recent results of \cite{CDS,T4} confirmed the  Yau-Tian-Donaldson conjecture
for smooth Fano manifolds. On the other hand, it will be interesting to investigate the K\"ahler-Einstein metrics on singular varieties. In \cite{DS1,DS2}, Donaldson-Sun showed that the Gromov-Hausdorff limit of a sequence of anti-canonically polarized  K\"ahler-Einstein manifold  is a $\mathbb Q$-Fano 
variety with Klt singularities and the tangent cone of any point of the limit variety is unique. Motivated by that, Li-Liu-Xu \cite{Li-Liu-Xu} formulated the stable degeneration conjecture for any algebraic varieties with Klt singularities. Also,  Hein-Sun \cite{H-S} proved that for each Calabi-Yau variety with isolated cone singularity, that is, the singularity is the cone over a smooth Fano Einstein vareity,  the global K\"ahler-Einstein metric is asymptotically the same as the local Ricci flat metric constructed by Calabi Ansatz.
This interesting result shows that among all local K\"ahler-Einstein metric, the one comes from Calabi Ansatz is most stable and canonical. Actually, the K-stability condition of the Fano Einstein manifolds is essentialy used in their proof. 
It is natural to ask similar questions for algebraic varietyies with log canonical singularities both globally and locally as follows: 

\begin{question}\label{q}Given a log-canonical cone singularity, can we construct a complete K\"ahler-Einstein metric locally as a model metric? If it exits, how canonical and stable it is compared to all the complete metrics induced from the global K\"ahler-Einstein metrics on compact algebraic varieties with such a singularity?
\end{question} 
We note that the technique of constructing model metrics with good curvature condition has been widely used in solving Monge-Amp\`ere equations in various geometric setting, for instance, Cheng-Yau's \cite{CY2} construction of complete K\"ahler-Einstein metric in strongly pseudoconvex domain,  Koabayashi's \cite{Kob} construction of complete K\"ahler-Einstein metric for an ample pair $K_X+D>0$ (see also \cite{TY},\cite{C-G}).

However, the construction of model metrics around singularities in Question $\ref{q}$ will be harder than the construction of model metrics for ample pairs in Koabayashi's setting. For if we take a log resolution of the singularity, the divisor $K_X+D$  is no longer ample but only semi-positive. A geometric explanation of this fact is that 
we do not only have the degeneration along the normal direction of divisor $D$, but also the degeneration along the tangential direction of $D$. In complex dimension two, the interesting work of Kobayashi \cite{Kob2} and Nakamura \cite{Nak} classified two dimensional log canonical singularity and showed that they are uniformized by bounded symmetric domains equipped with invariant Bergman metrics.  Recently Song \cite{S4} proves that for a family of canonically polarized varieties, the negative K\"ahler-Einstein metrics in the nearby fibers converge in the Gromov-Hausdorff limit to a singular K\"ahler-Einstein metric on the central fiber which has complete ends towards the locus of Non-Klt center. This indicates that it's interesting to study local complete K\"ahler-Einstein metrics with negative scalar curvature near non-Klt log canonical singularities. 
\medskip

In this short note, we aim to  construct  model metrics for  log-canonical cone singularities with a large symmetry in any dimension, which generalizes the work of Kobayashi \cite{Kob2} in complex dimension two.

%

\section{Singularities in algebraic geometry}

\begin{definition} \label{sing} Let $X$ be a normal projective variety such that $K_X$ is a $\mathbb{Q}$-Cartier divisor. Let $ \pi :  Y \rightarrow X$ be a log resolution and $\{E_i\}_{i=1}^p$ the irreducible components of the exceptional locus ${\rm Exc}(\pi)$ of $\pi$ with a unique collection $\{a_i\}\subset \mathbb{Q}$ such that
$$K_Y = \pi^* K_X + \sum_{i=1}^{ p } a_i E_i .$$ Then $X$ is said to have

\begin{enumerate}
\item[$\bullet$]  log terminal singularities if $a_i > -1$ for all $i $;
\medskip

\item[$\bullet$]  log canonical singularities if $a_i \geq -1$ for all $i$.
\medskip
\end{enumerate}
\end{definition}

\begin{example}\label{e1}
Let $E$ be an elliptic curve and $L$ a negative line bundle on $E$. Contracting the zero section of $L$, we get a singular variety $X$ with an isolated cone singularity. Then the contraction map $f:L\rightarrow X$ is a resolution of $X$ with zero section of $L$ the exception divisor which we still denote by $E$.
Write
$$K_Y+E=f^*K_X+aE$$
for some constant $a$. Restricting the above formula to $E$ and applying adjunction formula, we have $$K_E=0+aE|_{E}$$
Noticing that $K_E=0$, we get $a=0$. Therefore $X$ is log canonical by Definition \ref{sing}.  In fact the isolated singularity of $X$ is non-Klt log-canonical.
\end{example}

\begin{example}[See page $96$ in \cite{KK}]\label{e2}
In general, for a smooth Calabi-Yau manifold  $M$ and a negative line bundle $L$ on $M$, we can get a log canonical cone singular variety $\widehat M$ by contracting the zero section of $L$. Actually,
$$\widehat M={\rm Spec}_M \big(\mathbb \bigoplus_{i=0}^{i=\infty}H^0(M,-L^i)\big).$$
\end{example}

\section{The construction of complete K\"ahler-Einstein metrics near cone singularities}
\subsection{The approach via the invariant Bergman metric}
The uniformization of log canonical singularity by K\"ahler-Einstein metric in complex dimension two is due to Kobayashi \cite{Kob2}, our construction below is a high dimensional generalization. 
\begin{theorem}\label{t1}
Let  $A$ be an abeliean variety with complex dimension $n$ and $N$ a negative line bundle on $A.$ By contracting the zero section of $N,$  one obtains a singular variety $\widehat X$. 
Let   $(\widehat X,o)$ denote the germ of the isolated singularity of $\widehat X$.  Then there is an open neighborhood (in Euclidean topology) $\widehat {\mathcal U}$ of $o$ in $\widehat X$ such that ${\mathcal U}:=\widehat{\mathcal U}\backslash \{o\}$  is a smooth quotient space of a unit complex ball $\mathbb B^{n+1}$ by a discrete subgroup of ${\rm Aut}(\mathbb B^{n+1}).$ As a consequence, $\mathcal U$  has a negative K\"ahler- Einstein metric  induced from the Bergman metric of the ball which is complete towards $o$.
\end{theorem}
\noindent{\bf Proof: }
Denote by $f:N\rightarrow \widehat X$ the contraction map. For simplicity, we still denote the zero section of $N$ by $A$. Then $f^{-1}(o)=A.$

The idea of the proof is as follows. First, we will take the unit disc bundle over $A$, sitting inside $N$, with respect to a certain Hermitian metric $\rho$ on $N$; take $\widehat {\mathcal U}$ to be the image of the disc bundle under the contraction map $f$. Next, we will show that  the universal cover of $\widehat {\mathcal U}\backslash  A$ is a horoball of the  complex unit ball.   Furthermore, by constructing the quotient map explicitly, we  will prove that the deck transformations of the universal cover are induced by the biholomorphisms of the complex unit ball.  Therefore, the Bergman metric of the unit ball descends to a metric on $\widehat {\mathcal U}\backslash \mathcal A$ which is complete towards $A$. 

Now we will proceed to prove the theorem.

Following \cite{GH}, we first define a Hermitian metric $\rho$ on $N$ explicitly as follows. Let $A=\mathbb C^n/\Gamma:=V/\Gamma$ and $\pi:\mathbb C^n\rightarrow A$ is the covering map.  Noticing that $c_1(N)$ is an invariant integral form, negative of type (1,1), we choose a basis $\lambda_1,...,\lambda_{2n}$ for $\Gamma$ over $\mathbb Z$ such that in terms of dual coordinates $x_1,x_2,\cdots,x_{2n}$ on $V$
$$c_1(N)=\sum_{\alpha=1}^{n}\delta_{\alpha}dx_{\alpha}\wedge dx_{n+\alpha},\qquad\delta_{\alpha}\in\mathbb Z.$$
Since $c_1(N)$ is non-degenerate, we can further take as our basis for the complex vector space $V$ the vectors 
$$e_{\alpha}=\delta^{-1}_{\alpha}\lambda_{\alpha},\qquad\alpha=1,\cdots,n.$$
The period matrix of $\Gamma\subset V$ will then be of the form $\Omega=(\Delta_{\delta},Z)$, where $Z$ is a $n\times n$ matrix and 
\begin{displaymath}\Delta_{\delta}=\left(\begin{matrix}\delta_1&&&&0\\
&\cdot&&&\\
&&\cdot&&\\
&&&\cdot&\\
0&&&&\delta_n
\end{matrix}\right).
\end{displaymath}
We further write $c_1(N)$ in the matrix form as follows:
\begin{displaymath}Q=\left(\begin{matrix}0&\Delta_{\delta}\\
-\Delta_{\delta}&0
\end{matrix}\right).
\end{displaymath}
Noticing that $c_1(N)$ is a form of type $(1,1)$, we have $\Omega\cdot Q^{-1}\cdot \Omega^{t}=0,$ i.e,
\begin{displaymath}(\Delta_{\delta},Z)\left(\begin{matrix}0&\Delta_{\delta}\\
-\Delta_{\delta}&0
\end{matrix}\right)
\left(\begin{matrix}\Delta_{\delta}\\
Z^t
\end{matrix}\right)=Z-Z^t=0,
\end{displaymath}
which shows $Z$ is symmetric.  By the negativity of $c_1(N)$, we have $-\sqrt{-1}\cdot\Omega\cdot Q_{\delta}^{-1}\cdot\bar\Omega^t<0,$ i.e., 
\begin{displaymath}-\sqrt{-1}(\Delta_{\delta},Z)\left(\begin{matrix}0&-\Delta^{-1}_{\delta}\\
\Delta^{-1}_{\delta}&0
\end{matrix}\right)
\left(\begin{matrix}\Delta_{\delta}\\
\bar Z^t
\end{matrix}\right)=-\sqrt{-1}(Z-\bar Z^t)=2\Im Z<0.
\end{displaymath}
Therefore $Z=X-\sqrt{-1}Y$, where $X$ and $Y$ are $n\times n$ real symmetric matrices and $Y$ is positive definite. 
e
For later use, we  will write the matrix form of $c_1(N)$ in terms of  $dz_i\wedge \bar dz_j,i,j=1,\cdots,n$ . Here $z_i$ are the dual of $e_i$.    Notice that the period matrix is $\Omega=(\Delta_{\delta},Z)$ and 
$$(\lambda_1,\cdots,\lambda_{2n})^t=(e_1,\cdots,e_n)\Omega,$$
\begin{displaymath}(dz_1,\cdots,dz_n,d\bar z_1,\cdots,d\bar z_n)=
\left(\begin{matrix}\Omega\\
\bar\Omega
\end{matrix}\right)(dx_1,\cdots,dx_{2n})^t.
\end{displaymath}
By changing the basis from $\{dx_i\wedge dx_j\}$ to $\{dz_{\alpha}\wedge d\bar z_{\beta}\}$, we have 

$$c_1(N)=\frac{-\sqrt{-1}}{2}\sum_{\alpha,\beta=1}^{n}W_{\alpha\beta}dz_{\alpha}\wedge d\bar z_{\beta},$$
where $W=(W_{\alpha\beta})=Y^{-1}$. 



Since $\pi^*N$ is a trivial line bundle on $V$, we have a global trivialization 
$$\phi:\pi^*N\rightarrow V\times\mathbb C.$$  Now for $z\in V$ and $\lambda\in\Gamma, $ the fibers of $\pi^*N$ at $z$ and $z+\lambda$ are by definition both identified with the fiber of $N$ at $\pi(z),$ and comparing the trivialization $\phi$ at $z$ and $z+\lambda$ yields an automorphism of by $\mathbb C$:
$$\mathbb C\xleftarrow{\phi_z}(\pi^*N)_z=N_{\pi(z)}=(\pi^*N)_{z+\lambda}\xrightarrow{\phi_{z+\lambda}}\mathbb C.$$
Such an automorphism is given as multiplication by a nonzero complex number. If we denote this number by $e_{\lambda}(z),$ we obtain a collection of transition functions $\{e_{\lambda}\in\mathcal O^*(V)\}_{\lambda\in\Gamma}$.  The functions $e_{\lambda}$ necessarily satisfy the compatibility relation
$$e_{\lambda^{\prime}}(z+\lambda)e_{\lambda}(z)=e_{\lambda}(z+\lambda^{\prime})e_{\lambda^{\prime}}(z)=e_{\lambda+\lambda^{\prime}}(z)$$
for all $\lambda,\lambda^{\prime}\in\Gamma.$

We choose  a collection of transition functions associated to the line bundle $N$ as follows:
$$e_{\lambda_{\alpha}}(z)\equiv 1,e_{\lambda_{n+\alpha}}(z)\equiv e^{-2\pi i (z_{\alpha}-t_{\alpha})},\alpha=1,\cdots,n,$$
where $z_{\alpha}$ is the $\alpha$-th component of $z$ and $t_{\alpha}$ is a constant. (Here $t_{\alpha}$ is determined by the line bundle.) For any local section $\tilde\theta$ of $N$ over $U\subset A$, $\theta=\phi^*(\pi^*\tilde\theta)$ is an analytic function on $\pi^{-1}(U)$ satisfying 
$$\theta(z+\lambda_{\alpha})=\theta(z),\quad\theta(z+\lambda_{n+\alpha})=e^{-2\pi i(z_{\alpha}-t_{\alpha})}\theta(z).$$  Then define the Hermitian metric $||\,\,||$ on $N$ as 
$$||\tilde\theta(z)||^2=h(z)|\theta(z)|^2, \quad {\rm where}\,\, h(z)=e^{(-\pi/2)\sum W_{\alpha\beta}(z_{\alpha}-\bar z_{\alpha}+iY_{\alpha\alpha}-t_{\alpha}+\bar t_{\alpha})(z_{\beta}-\bar z_{\beta}+iY_{\beta\beta}-t_{\beta}+\bar t_{\beta})}.$$
It is easy to verify that $||\tilde\theta(z+\lambda)||=||\tilde\theta(z)||$ for all $\lambda\in \Gamma$.  Hence $||\,\,||$ is a well defined metric on $N$, which we denote by $\rho$.

Define $\mathcal D_k:=\{w\in N;0<\rho(w,w)<e^{-k}\}$,  $\mathcal  W_k:=\{(w,z)\in\mathbb C\times\mathbb C^n\,|0<|w|^2h(z)<e^{-k}\}$ and  $\mathcal H_k:=\{(u,z)\in\mathbb C\times\mathbb C^n\,|-\Im u-\frac{\pi}{2}\sum W_{\alpha\beta}(z_{\alpha}-\bar z_{\alpha}+iY_{\alpha\alpha}-t_{\alpha}+\bar t_{\alpha})(z_{\beta}-\bar z_{\beta}+iY_{\beta\beta}-t_{\beta}+\bar t_{\beta})<-k\}$.  We have the following commutative diagram:

\begin{equation}
\begin{array}{cccccc}
&\mathcal H_k& \xrightarrow{\exp_k}& \mathcal W_k &\xrightarrow {\pi_k}& \mathcal D_k\\
&\Big\downarrow\rlap{$\scriptstyle i_1$}&&\Big\downarrow\rlap{$\scriptstyle i_2$}&&\Big\downarrow\rlap{$\scriptstyle i_3$}\\
&\mathbb C^{n+1}& \xrightarrow{\exp}& \mathbb C^{n+1} &\xrightarrow{\pi}& N \\
\end{array}.
\end{equation}
Here $i_1,i_2,i_3$ are natural injections; ${\rm exp:\mathbb C^{n+1}\rightarrow\mathbb C^{n+1}}$ is defined by $(u,z)\mapsto(e^{\frac{i}{2}u},z)$; ${\rm exp}_k$ and $\pi_k$ are the restrictions of ${\rm exp}$ and $\pi$ respectively.  We will show  $\mathcal H_k$ is biholomorphic to a horoball, hence it is the universal cover of $D_k$.  Notice that the deck transformations group of  $\mathcal H_k \xrightarrow{\pi_k\circ \exp_k}\mathcal D_k $ can be presented as a subgroup of the matrix group as follows:
\begin{displaymath}\left\{
\left(\begin{matrix}1&\quad\quad -4\pi  m_1&\cdots&\quad\quad -4\pi  m_n & 2\pi\sum m_{\alpha}Z_{\alpha\alpha}-2\pi \sum Z_{\alpha\beta}m_{\alpha}m_{\beta}+4p\pi \\
0&1&\cdots&0&\delta_1l_1+Z_{1\alpha}m_{\alpha}\\
\cdots&\cdots&\cdots&\cdots&\cdots\\
0&0&1&\cdots&\delta_nl_{n}+Z_{n\alpha}m_{\alpha}\\
0&0&\cdots&0&1
\end{matrix}
\right);m_{\alpha},l_{\alpha},p\in\mathbb Z\right\}.
\end{displaymath}
Here $\Omega=(\Delta_{\delta},Z)$; $Z_{\alpha \beta}$ and $\delta_{\alpha}$ are entries of $\Delta_{\delta}$ and $Z $ respectively.
Moreover, the action of an element of the group on $\mathbb C^{n+1}$ can be presented as matrix multiplications as follows:
\begin{displaymath}
\left[\begin{matrix}u^{\prime}\\z_1^{\prime}-t_{1}\\\cdot\\\cdot\\z_{n}^{\prime}-t_{n}\\1\end{matrix}\right]=
\left(\begin{matrix}1&\quad\quad -4\pi  m_1&\cdots&\quad\quad -4\pi  m_n & 2\pi\sum m_{\alpha}Z_{\alpha\alpha} -2\pi \sum Z_{\alpha\beta}m_{\alpha}m_{\beta}+4k\pi \\
0&1&\cdots&0&\delta_1l_1+\sum Z_{1\alpha}m_{\alpha}\\
\cdots&\cdots&\cdots&\cdots&\cdots\\
0&0&1&\cdots&\delta_nl_{n}+\sum Z_{n\alpha}m_{\alpha}\\
0&0&\cdots&0&1
\end{matrix}
\right)\left[\begin{matrix}u\\z_1-t_{1}\\\cdot\\\cdot\\z_{n}-t_{n}\\1\end{matrix}\right]
\end{displaymath}

Notice that this action preserves $|w|^2h(z)$, or equivalently, leaves invariant the following function:
$$-\Im u+(-\pi/2)\sum W_{\alpha\beta}(z_{\alpha}-\bar z_{\alpha}+iY_{\alpha\alpha}-t_{\alpha}+\bar t_{\alpha})(z_{\beta}-\bar z_{\beta}+iY_{\beta\beta}-t_{\beta}+\bar t_{\beta});$$
Hence this action extends to a biholomorphism for each $\mathcal H_k,k\in\mathbb R$.

On the other hand, since $W$ is positive definite, there is an $n\times n $ matrix $V$ such that $W=V\bar V^t.$  Then we can perform coordinates change of $\mathbb C^{n+1}$  as follows:
\begin{equation}\label{coor}
\begin{split}
&\tilde u=u+i\pi W_{\alpha\beta}z_{\alpha}z_{\beta}+\frac{\pi}{2}W_{\alpha\beta}(iY_{\alpha\alpha}-t_{\alpha}+\bar t_{\alpha})(iY_{\beta\beta}-t_{\beta}+\bar t_{\beta})+2\pi i\sum W_{\alpha\beta}z_{\alpha}(iY_{\beta\beta}-t_{\beta}+\bar t_{\beta}),\\
&\tilde z_{\alpha}=\sum V_{\alpha\beta}z_{\beta},\,\,\alpha=1,2,\cdots,n.
\end{split}
\end{equation}
Under this coordinates change, the domain $\mathcal H_k=\{(u,z)|-\Im u+(-\pi/2)\sum W_{\alpha\beta}(z_{\alpha}-\bar z_{\alpha}+iY_{\alpha\alpha}-t_{\alpha}+\bar t_{\alpha})(z_{\beta}-\bar z_{\beta}+iY_{\beta\beta}-t_{\beta}+\bar t_{\beta})<-k\}$ is mapped to the domain
$\{(\tilde u,\tilde z)|\Im {\tilde u}-\sum_{i=1}^n|\tilde z|^2>k\}$ biholomorphically.  Particularly, when $ k=0$, the latter one is the Heisenberg model of the unit ball in $\mathbb C^{n+1}.$  

Since Bergman metric is invariant under automorphism group actions, the Bergman metric on the ball $\mathcal H_0$ descends to $\mathcal D_0$ and  its restrictions on $\mathcal H_k,k>0$ descend to $\mathcal D_k$.   Moreover, since the zero section of $N$ corresponds to the infinite boundary point of $\mathcal H_0$, that is, $\Im u=+\infty$, we further conclude that the metric on $\mathcal D_k$ is complete towards $A$.  Take $\mathcal U$ to be $\mathcal D_k$ for certain $k>0$, we conclude Theorem \ref{t1}.\,\, $\endpf$
\medskip

For the future investigation, we will construct a system of quasi-coordinates of $(\mathcal U,g)$ in the following, where $g$ is the metric constructed in Theorem \ref{t1}. 
Firstly, we recall the following definition of quasi-coordinates according to \cite{CY2} and \cite{Kob}.
\begin{definition} Let $V$ be a domain in $\mathbb C^m$. Let $X$ be an $m$-dimensional complex manifold and $\phi$ a holomorphic map of $V$ into $X$. $\phi$ is called a quasi-coordinate map if $\phi$ is of maximal rank everywhere. In this case,
$(V,\phi)$ is called a quasi-coordinate of $X.$
\end{definition}
Adapting to our case, we further make the following definition.
\begin{definition}\label{quasi}Let $\mathcal U=\widehat {\mathcal U}\backslash\{o\}$ as above. Assume $g$ is a complete K\"ahler metric on $\mathcal U$ towards $o$.  Then a system of quasi-coordinates of $(\mathcal U,g)$ is a set of quasi-coordinates $\Gamma=\{(V_{\alpha},\phi_{\alpha},(v_{\alpha}^1,v_{\alpha}^2\cdots, v_{\alpha}^{n+1}))\}$ of $\mathcal U$ with the following properties:
\begin{enumerate}[label=(\alph*)]
\item $\mathcal U=\bigcup_{\alpha}({\rm Image} \,\,{\rm of} \,V_{\alpha})$;
\item The complement of a certain open neighborhood $U\subset\mathcal U$ of the infinity point $o$ is covered by a finite number of quasi coordinates which are coordinate charts in the usual sense;
\item For each point $x\in U$, there is a quasi-coordinate $V_{\beta}$ and $\tilde x\in V_{\beta}$, such that $\phi_{\beta}(\tilde x)=x$ and dist$(\tilde x,\partial V_{\beta})\ge \epsilon_1$ in the Euclidean sense, where $\epsilon_1$ is constant independent of $\beta$;
\item There are positive constant $c$ and $A_k,k=1,2,\cdots,$ independent of ${\alpha}$, such that for each quasi coordinate $(V_{\alpha},\phi_{\alpha},(v_{\alpha}^1,v_{\alpha}^2\cdots, v_{\alpha}^{n+1}))$,  the following inequalities hold:
$$c^{-1}(\delta_{i\bar j})\le (g_{\alpha i\bar j})\le c(\delta_{i\bar j})$$
$$ |\frac{\partial^{p+q}}{\partial v_{\alpha}^p\bar v_{\alpha}^q}g_{\alpha i\bar j}|<A_{p+q},\forall p,q,$$ 
where $(g_{\alpha i\bar j})$ denote the metric tensor with respect to $(V_{\alpha},\phi_{\alpha},(v_{\alpha}^1,v_{\alpha}^2\cdots, v_{\alpha}^{n+1}))$.
\end{enumerate}
\end{definition}
\begin{lemma}\label{qc}  $(\mathcal U,g)$  has a system of quasi-coordinates.
\end{lemma}	
\noindent{\bf Proof:} Without loss of generality, we can assume that $\mathcal U=\mathcal D_1$.  Take for $\mathcal H_0$ the coordinates $(\tilde u,\tilde z_1,\cdots,\tilde z_n)$ as constructed in formula $(\ref{coor})$. 

For the later use, we introduce some special biholomorphisms. Denote by $\iota$ the following biholomorphic map
\begin{equation*}
\begin{split}
&\mathcal H_k\rightarrow\mathcal H_k\\
&\tilde u\mapsto \tilde u+4\pi,\tilde z\mapsto \tilde z\,\,.
\end{split}
\end{equation*}
Denote by  $\tau_a$,  $a>0$, the biholomorphic map of $\mathcal H_0$ as follows:
\begin{equation*}
\begin{split}
&\mathcal H_0\rightarrow\mathcal H_0\\
&\tau_a:\tilde u\mapsto a^{-2}\tilde u,\tilde z\mapsto a^{-1}\tilde z\,\,.
\end{split}
\end{equation*}

Take a fundamental domain $F$ of $V/\Gamma$ such that it is a bounded convex polytope containing the origin of $V$. Define the following bounded domains contained  in $\mathcal H_0\subset\mathbb C^{n+1}$ as
\vspace{-0.1in}
$$T:=\{(\tilde u,\tilde z)|1<\Im {\tilde u}-\sum_{i=1}^n|\tilde z|^2<4,0<\Re \tilde u<8\pi, \tilde z\in 3F \},$$
\vspace{-0.1in}
$$T':=\{(\tilde u,\tilde z)|2<\Im {\tilde u}-\sum_{i=1}^n|\tilde z|^2<3,0<\Re \tilde u<4\pi, \tilde z\in 2F \}.$$
Here we denote by $mF$, $m>0$, the set $\{(ma_1,\cdots,ma_n)\big|(a_1,\cdots,a_n)\in F\}$.
It is easy to see that $T'\subset T$  and dist$(\partial T,\partial T')>\epsilon>0$. 
\medskip

{\bf Claim:} For each point $q \in \mathcal D_2$, there is a point $Q\in T'$, an integer $m$, a real number $a>1$ and an element $\gamma$ of $\Gamma$ such that $(\exp_1\circ\pi_1\circ\tau_a\circ\gamma\circ l^m)(Q)=q.$

{\bf Proof of the claim:} Take a preimage $Q$ of $q$ in $\mathcal H_2$. Then there is an integer $m$ such that $0<\Re (\tilde u(l^m(Q)))<4\pi$, where $\tilde u(l^m(Q))$ is the $\tilde u$ coordinates of the point $l^m(Q)$. Next, there is an element $\gamma$ of the lattice $\Gamma$ such that the $\tilde z$ coordinates of $(\gamma\circ l^m)(Q)$ is contained in $2F$. At last, there is a real number $a>1$ such that $2<\Im {\tilde u((\tau_a\circ \gamma\circ l^m)(Q))}-\sum_{i=1}^n|\tilde z((\tau_a\circ\gamma\circ l^m)(Q))|^2<3$. Since $F$ is convex and $a>1$, the $\tilde z$ coordinates of $(\tau_a\circ\gamma\circ l^m)(Q)$ is still contained in $2F$.  We complete the proof of the claim.\,\,$\endpf$
\medskip

Next, we define our system of quasi-coordinates of $(\mathcal U,g)$as follows:
\begin{enumerate}
	\item for each point $q\in \mathcal D_2$, take $V_{q}$ to be $T^{\prime}\subset\mathbb C^{n+1}$;
	\item take coordinates $(v_{q}^1,v_{q}^2\cdots, v_{q}^{n+1})$ to be the coordinates $(\tilde z_1,\cdots,\tilde z_n, \tilde u)$ of $\mathcal H_1\subset\mathbb C^{n+1}$;
	\item take $\phi_{q}$ to be the restriction of the map $\exp_k\circ\pi_k\circ l^{-m_q}\circ\gamma_q^{-1}\circ\tau_{a_q}^{-1}$ to $T^{\prime}$, where the maps $\tau_{a_q},\gamma_q$ and $l^{m_q}$ are the those constructed in the claim associated with $q$.
\end{enumerate} 
It is clear that properties $(a),(b)$ in Definition \ref{quasi} is satisfied by our system. Because the curvature of Bergman metric is bounded, the curvature property $(d)$ holds as well. Take $U$ to be $\mathcal D_2$, it is clear that the property $(c)$ is also satisfied by our system.  Therefore, we complete the proof of Lemma $\ref{qc}$ .\,\,$\endpf$

\begin{remark} If we compactify the example in (\ref{e2}), we get a compact singular variety $\widehat X^{'}$ with an isolated cone singularity. According to \cite{Kob,Ye}, we can construct a global complete K\"ahler-Einstein metric on  $\widehat X^{'}$ by solving the corresponding
Monge-Amp\`ere equation.

\end{remark}

\subsection{The approach via Calabi Ansatz}
\begin{theorem}\label{t2}
Let  $M$ be a Calabi-Yau manifold  with complex dimension $n$ and $\widehat M$  the total space of a negative line bundle $L$ on $M.$ By contracting the zero section of $L,$  one obtains a singular variety $\widehat X$. 
Let   $(\widehat X,o)$ denote the germ of the isolated singularity of $\widehat X$.  Then there is an open neighborhood (in Euclidean topology) $\widehat {\mathcal U}$ of $o$ in $\widehat X$ such that ${\mathcal U}:=\widehat{\mathcal U}\backslash \{o\}$ carrys a negative complete K\"ahler- Einstein metric.

\end{theorem}
\noindent{\bf Proof:}
Since $C_1(L)<0$, by Yau's theorem, there is a Ricci flat metric $\omega_M$ in class $-C_1(L).$ Choose a Hermitian metric $h$ on $L$ such that $Ric_h=-\omega_M$. Let $r$ be the distance function on the line bundle; locally $r=(a(z)|\eta|^2)$ where $z=(z_1,\cdots,z_n)$ are coordinates on $M$, $\eta$ is the coordinate on $L$ and $a(z)$ is the negative hermitian metric on $L$.

Next we will proceed by using a Calabi Ansatz type argument. Suppose that the potential metric can be written as $\omega=\omega_M+\ddbar\rho(r)$, where $\rho$ is a one variable function to be determined.

In local coordinates $\omega$ is expressed as
$$\omega=(1+\rho_rr)\omega_M+a(\rho_r+\rho_{rr}a\bar \eta \eta)\nabla \eta\wedge\bar {\nabla \eta},$$
where $\nabla \eta=d\eta+a^{-1}\partial a\eta$.
Note that ${dz^i,\nabla\eta}$ are dual to the basis consisting of horizontal and vertical vectors
$$\nabla_{z^i}=\frac{\partial}{\partial z^i}-a^{-1}\frac{\partial a}{\partial z^i}\eta\frac{\partial}{\partial\eta}\,\,{\rm and}\,\,\frac{\partial}{\partial\eta}.$$

Notice that $\omega$ is positive if and only if 
$$1+\rho_rr>0,\rho_r>0,\rho_{rr}r+\rho_r>0.$$ 
Calculation yields that
$$\det\omega=(1+r\rho(r))^n\omega_M^na(\rho_{rr}r+\rho_r)d\eta\wedge d\bar\eta.$$
Hence the K\"ahler-Einstein equation $Ric_{\omega}=-\omega$ reads as
\begin{equation*}
\begin{split}
-\ddbar\log\det\omega&=-\ddbar\log(1+r\rho)^n(\rho_{rr}r+\rho_r)-\ddbar\log\omega_M^n-\ddbar\log a\\
&=-\omega_M-\ddbar\rho.
\end{split}
\end{equation*}
Simplifying it,  we get 
$$\ddbar\log(1+r\rho)^n(\rho_{rr}r+\rho_r)=\ddbar\rho.$$
In order to derive a solution, it suffices to solve the following ODE \begin{equation}\label{ODE}(1+r\rho)^n(\rho_{rr}r+\rho_r)=e^{\rho}.
\end{equation}
Now let us make a change of variables by setting $s=\log r$. Noticing the fact that $$r\rho_r=\rho_s\,\,{\rm and}\,\,(\rho_{rr}r+\rho_r)=(\rho_rr)_r=(\rho_s)_r=\frac{\rho_{ss}}{e^s},$$
we have that the equation (\ref{ODE}) is equivalent to 
$$(1+\rho_s)^n\rho_{ss}=e^{\rho+s}.$$
Differentiating both sides, we get
$$n(1+\rho_s)^{n-1}\rho_{ss}^2+(1+\rho_s)^n\rho_{sss}=e^{\rho+s}(\rho_s+1)=\rho_{ss}(1+\rho_s)^{n+1}.$$
Defining $f:=\rho_s+1$, we have $$\frac{nf_s}{f}+\frac{f_{ss}}{f_s}=f.$$
Since $f_s>0$, we can solve $s$ as function of $f$ by inverse function theorem. Letting $s=\Phi(f)$ and denoting $f_s(\Phi(f))$ by $g(f)$ , we transform the above ODE into the following form:
$$\frac{ng}{f}+g_f=f,$$ which implies that
$$(gf^n)'=f^{n+1}.$$
Solving it, we have $$f_s=g(f)=\frac{f^2}{n+2}+\frac{C}{f^n}$$
Letting $C=0$, we have $$\rho_s=f+1=-\frac{n+2}{s}-1.$$
Therefore  $\rho=-(n+2)\ln |s|-s$ is a solution.
Since we are consider a neighborhood of the vertex, $s=\log r$ is negative here.  As a conclusion, we derive the desired metric $$\omega=\omega_M+\ddbar(\log\log|r|^{-(n+2)}-\log r)=\ddbar\log\log|r|^{-(n+2)}.\,\,\,\,\,\endpf$$

\begin{remark}
The following computation shows that if the Calabi-Yau metric $\omega_M$ in class $-C_1(L)$ is not flat, the curvature of metric $\omega$ will blow up. Firstly we have that 
$$\omega_M=\ddbar\log a=\frac{a_{z\bar z}}{a}-\frac{|a_z|^2}{a^2}dz\wedge d \bar z.$$
Calculation yields that
\begin{equation*}
\begin{split}
\omega=(\frac{-3}{\log r}(\frac{a_{z\bar z}}{a}&-\frac{|a_z|^2}{a^2})+\frac{3}{(\log r)^2}\frac{|a_z|^2}{a^2})dz\wedge d\bar z+\frac{3a_{\bar z}\bar\eta}{r(\log r)^2} d\eta\wedge d\bar z\\
&+\frac{3a_z\eta}{r (\log r)^2} dz\wedge d\bar \eta+\frac{3a}{r(\log r)^2}d\eta\wedge d\bar\eta.
\end{split}
\end{equation*}
Notice that in the $z$-direction the leading term of the metric is  $\frac{-3}{\log r}(\frac{a_{z\bar z}}{a}-\frac{|a_z|^2}{a^2})dz\wedge d\bar z$, hence in most cases (when $\omega_M$ is not very special) the curvature will blow up in the rate of $|\log r|$ if $(\frac{a_{z\bar z}}{a}-\frac{|a_z|^2}{a^2})$ is not a very simple function.
\end{remark}

\begin{remark}
The metric $\omega$ constructed above seems similar to the Carlson-Griffiths metric used for instance in \cite{C-G}, but there are two differences: Firstly, our metric degenerates also along the tangential direction of $M$ while Carlson-Griffiths doesn't; secondly, our potential $\log|r|^{-(n+2)}$ 
itself is a plurisubharmonic function locally while the one in Carlson-Griffiths is not.
\end{remark}

At last, we prove the following proposition.
\begin{proposition} The metric constructed  in Theorem \ref{t1} coincides the metric constructed in  Theorem \ref{t2}.
\end{proposition}
\noindent{\bf Proof:}
  Note that the metric in Theorem \ref{t1} is actually
  $$\ddbar\log|\log(|w|^2h(z))|,$$ where $z=(z_1,\cdots,z_n)$ are the coordinates on the abeliean variety and $w$ is the coordinate on the line bundle. Hence $|w|^2h(z)$ is the distance function. On the other hand, the metric in Theorem \ref{t2} takes the form of 
  $$\ddbar\log\log|r|^{-(n+2)},$$ where $r$ is the distance function. This finishes the proof.\,\,$\endpf$

\section{Acknowledgement}
 We want to thank our advisors professor Xiaojun Huang and Jian Song for continuous encouragement. We also thank professor Martin de Borbon for notifying us the Calabi Ansatz method has been used in his recent interesting joint work with Professor Spotti in \cite{BS}.
 The second author also wants to thank professor Hans-Joachim Hein for teaching him a lot about the geometry of complete Calabi-Yau manifolds. 



\end{document}